\documentclass[10pt]{amsart}

\usepackage{geometry}
\usepackage{graphicx} % see geometry.pdf on how to lay out the page. There's lots.
\usepackage{datetime}
\usepackage{bookmark}
\usepackage{amsmath,amssymb,amsthm,mathrsfs,pdfpages}  % Better maths support & more symbols
\usepackage{pdfsync}  % enable tex source and pdf output syncr
\newtheorem*{theorem*}{Theorem}
\newtheorem{theorem}{Theorem}
\newtheorem{lemma}{Lemma}
\newtheorem{proposition}{Proposition}
\newtheorem*{corollary*}{Corollary}
\newtheorem{corollary}{Corollary}
\newtheorem{definition}{Definition}

\theoremstyle{remark}
\newtheorem*{remark*}{Remark}
\newtheorem{remark}{Remark}

\newtheorem*{convention*}{Conventions}

\newcommand{\thm}[1]{\begin{theorem}#1\end{theorem}}

\newcommand{\lem}[1]{\begin{lemma}#1\end{lemma}}

\newcommand{\eq}[1]{\begin{equation}#1\end{equation}}
\newcommand{\eqs}[1]{\begin{equation*}#1\end{equation*}}
\newcommand{\alig}[1]{\begin{align}#1\end{align}}
\newcommand{\aligs}[1]{\begin{align*}#1\end{align*}}
\newcommand{\nn}{\nonumber}
\newcommand{\pf}[1]{\begin{proof}#1\end{proof}}

\newcommand{\mbb}[1]{\mathbb{#1}}
\newcommand{\mcal}[1]{\mathcal{#1}}
\newcommand{\vGamma}{\mathit{\Gamma}}

\newcommand{\mmod}{\mathrm{mod}}
\newcommand{\ts}{\textstyle}
\newcommand{\ds}{\displaystyle}
\newcommand{\sss}{\scriptscriptstyle}
\newcommand{\dd}{\hspace*{0.1em}\mathrm{d}}% rectified d in integrals
\newcommand{\f}[2]{\frac{#1}{#2}}% Fractions
\newcommand{\el}{\left}% extending left brackets
\newcommand{\er}{\right}% extending right brackets
\newcommand{\re}{\mathrm{Re}}% Real part Re
\newcommand{\res}[1]{\underset{#1}{\,\,\mathrm{res}\,\,}}% Residue
\newcommand{\rarrow}{\rightarrow}

\newcommand{\lbar}[1]{\overline{#1}}
\newcommand{\sym}{\mathrm{sym}}
\newcommand{\sumprime}[1]{\sideset{}{'}\sum_{#1}}
\newcommand{\ep}{\varepsilon}

\newcommand{\hf}{\f{1}{2}}
\newcommand{\thf}{{\textstyle\f{1}{2}}}

\newcommand{\tf}[2]{{\textstyle\frac{#1}{#2}}}

\newcommand{\heq}{\hspace*{3.1ex}}

\newcommand{\qand}{\quad\mbox{and}\quad}
\newcommand{\subs}[1]{\substack{#1}}

\newcommand{\ints}[1]{\int_{\sss #1}}
\newcommand{\chiD}{\chi_{\sss D}}

\begin{document}

%%%%% To ease editing, for IMPAN journals add:

\baselineskip=17pt
%%%%%%%%%%%

%% In the running head, replace first names by initials
%% and give an abbreviation of the title.
\title[Symmetric Square $L$-functions]{The first moment of central values of symmetric square $L$-functions in the weight aspect}

\author[S. Liu]{Shenhui Liu}
\address{Department of Mathematics, The Ohio State University\\ 231 W 18th Avenue\\
Columbus, Ohio 43210, United States of America}
\email{liu.2076@osu.edu}

%\date{September 29, 2016}

%\newsavebox{\discard}
%\sbox{\discard}{\vbox{\tableofcontents}}

\begin{abstract}
    In this note we investigate the behavior at the central point of the symmetric square $L$-functions, the most frequently used $\rm{GL}(3)$ $L$-functions. We establish an asymptotic formula with arbitrary power saving for the first moment of $L(\hf,\sym^2f)$ for $f\in\mcal{H}_k$ as even $k\rarrow\infty$, where $\mcal{H}_k$ is an orthogonal basis of weight-$k$ Hecke eigencuspforms for $SL(2,\mbb{Z})$. The approach taken in this note allows us to extract two secondary main terms from the error term $O(k^{-\hf})$ in previous studies. More interestingly, our result exhibits a connection between the symmetric square $L$-functions and quadratic fields, which is the main theme of Zagier's work \textit{Modular forms whose coefficients involve zeta-functions of quadratic fields} in 1977. Specifically, the secondary main terms in our asymptotic formula involve central values of Dirichlet $L$-functions of characters $\chi_{-4}$ and $\chi_{-3}$ and depend on the values of $k\,(\mmod\ 4)$ and $k\,(\mmod\ 6)$, respectively.
\end{abstract}



\maketitle
%\tableofcontents
\section{Introduction}
Let $k>0$ be a large even integer and consider the Hilbert space $S_k$ of holomorphic cusp forms of weight $k$ for $\vGamma_0(1)=SL(2,\mbb{Z})$, with respect to the Petersson inner product $(\cdot,\cdot)$ on $\vGamma_0(1)\backslash\mbb{H}$. Fix a normalized orthogonal Hecke basis $\mcal{H}_k$ of $S_k$, which is of size $\#\mcal{H}_k=\f{1}{12}(k-1)+O(1)$. Here the normalization means that the first Fourier coefficient of forms in $\mcal{H}_k$ is $1$. In this note we are interested in the frequently used symmetric square $L$-function for $f\in\mcal{H}_k$
$$
L(s,\sym^2f)=\zeta(2s)\sum_{n\geq1}\f{\lambda_f(n^2)}{n^s}\quad(\re(s)>1),
$$
which is a $\rm{GL}(3)$ $L$-function by Gelbart$-$Jacquet \cite{GelbartJacquet1978}. By Shimura \cite{Shimura1975}, Zagier \cite{Zagier1977}, and more generally \cite{GelbartJacquet1978}, $L(s,\sym^2f)$ has analytic continuation to the whole $s$-plane.

Several authors (Lau \cite{Lau2002}, Khan \cite{Khan2007}, and Sun \cite{Sun2013}) obtained asymptotic formulas for the first moment of central value of symmetric square $L$-functions by using Petersson's formula and various analytic techniques, among which Sun's result is the strongest:
\aligs{
\sum_{f\in\mcal{H}_k}{\rm{w}}_f L\Big(\hf,\sym^2 f\Big)=\psi\Big(k-\hf\Big)+2\gamma+\hf\psi\Big(\f{3}{4}\Big)-\log(2\pi^{\f{3}{2}})+O(k^{-\hf}),
}
where
$${\rm{w}}_f=\f{\Gamma(k-1)}{(4\pi)^{k-1}(f,f)},$$
$\psi(z)=\f{\Gamma'(z)}{\Gamma(z)}$ denotes the digamma funciton and $\gamma$ is the Euler constant. Still based on application of Petersson's formula, we use a different approach to establish an improved asymptotic formula with arbitrary power saving.

\thm{\label{MainThm}
For any $B>0$ and sufficiently large even integer $k>0$ we have
\alig{
\nn&\sum_{f\in\mcal{H}_k}{\rm{w}}_f L\Big(\hf,\sym^2 f\Big)=\\
\tag{${\rm{M}}_{1}$}&\heq\psi\Big(k-\hf\Big)+2\gamma+\hf\psi\Big(\f{3}{4}\Big)-\log(2\pi^{\f{3}{2}})\\
\tag{${\rm{M}}_{-4}$}&\heq+i^{-k}\sqrt{\f{\pi}{2}}L\Big(\hf,\chi_{-4}\Big)\f{\Gamma\!\el(\f{k-\hf}{2}\er)}{\Gamma\!\el(\f{k+\hf}{2}\er)}\\
\tag{${\rm{M}}_{-3}$}&\heq+\sqrt{2\pi}i^{-k}L\Big(\hf,\chi_{-3}\Big)
\Big(\f{2}{\sqrt{3}}\Big)^{k-\hf}F\bigg(\f{k-\hf}{2},\f{k-\hf}{2};\hf;-\f{1}{3}\bigg)
\f{\Gamma\!\el(\f{k-\hf}{2}\er)}{\Gamma\!\el(\f{k+\hf}{2}\er)}\\
\nn&\heq+O_{\sss B}(k^{-B}).
}
Here $\chiD({\boldmath\cdot})=(\f{D}{{\ts\boldmath\cdot}})$ denotes the Kronecker symbol and $F(a,b;c;z)$ is the Gauss hypergeometric function.
}
In the process of preparing this work for publication, the author learned that Balkanova and Frolenkov \cite{BalkanovaFrolenkov2016} independently obtained a slightly stronger result by using a method different from ours. \\

Now we give the plan of the rest of the note with some remarks. Our approach is different from the above-mentioned studies in that we use a different approximate functional equation (Lemma \ref{AFE_SymmSq}) at the very beginning of all analysis. After including necessary ingredients and tools in \S\,2 we obtain the primary main term ${\rm{M}}_{1}$ in \S\,3. The new approximate functional equation makes a real difference in \S\,4 where we analyze the off-diagonal contribution: it gives us a handle to reduce the analysis at various stages to the problem of counting solutions of certain quadratic congruence equations, and thus allows us to extract the secondary main terms ${\rm{M}}_{-4}$ and ${\rm{M}}_{-3}$. The labeling of the main terms ${\rm{M}}_{D}$ with $D$ a discriminant is used to indicate that the corresponding analysis involves the zeta function $\zeta_{K_{\sss D}}(s)$ of the quadratic field $K_{D}$ with discriminant $D$, where we treat $\mbb{Q}$ as the degenerate ``quadratic'' field with discriminant $1$ (see (\ref{L_zeta}), Lemma \ref{L_N}, and Lemma \ref{L_M}). We also note that the central values of Dirichlet $L$-functions in ${\rm{M}}_{-4}$ and ${\rm{M}}_{-3}$ are nonzero. In fact one can use $Mathematica^\circledR$ to see
$$
L\Big(\hf,\chi_{-4}\Big)\approx 0.667691\qand L\Big(\hf,\chi_{-3}\Big)\approx0.480868.
$$
We remark that the treatment of the off-diagonal contribution is inspired by Lau$-$Tsang \cite{LauTsang2005}, in which arbitrary power saving is obtained in a different context.

Now we need to say more about ${\rm{M}}_{-4}$ and ${\rm{M}}_{-3}$. The main term ${\rm{M}}_{-3}$ looks bizarre at first glance but can be shown by further analysis (\S\,4.3.1) to be
\alig{
\tag{$\rm{M}^{'}_{-3}$}3^{\f{1}{4}}\sqrt{2\pi}L\Big(\hf,\chi_{-3}\Big)\f{\Gamma(k-\hf)}{\Gamma(k)}\big[S(k)+O(k^{-1})\big],
}
where
$$
S(k)=\begin{cases}
-1,&\quad\mbox{if }k\equiv2\,(\mmod\ 6),\vspace*{1ex}\\
0,&\quad\mbox{if }k\equiv4\,(\mmod\ 6),\vspace*{1ex}\\
1,&\quad\mbox{if }k\equiv0\,(\mmod\ 6).
\end{cases}
$$
By Barnes' formula \cite[(1.18.12)]{Erdelyi1953} %that for any fixed $a$ and large $z$ not in the vicinity of the negative real axis
%$$
%\log\Gamma(z+a)=\Big(z+a-\hf\Big)\log z-z+\hf\log2\pi+O(|z|^{-1}),
%$$
we have for large $k$
\aligs{
\f{\Gamma\!\el(\f{k-\hf}{2}\er)}{\Gamma\!\el(\f{k+\hf}{2}\er)}=\sqrt{2}k^{-\hf}e^{O(k^{-1})}\qand \f{\Gamma(k-\hf)}{\Gamma(k)}=k^{-\hf}e^{O(k^{-1})}.
}
Hence the secondary main term ${\rm{M}}_{-4}$ is of order exactly $k^{-\hf}$, and so is ${\rm{M}}_{-3}$ when $k\not\equiv4\,(\mmod\ 6)$; also ${\rm{M}}_{-3}$ is $O(k^{-\f{3}{2}})$ if $k\equiv4\,(\mmod\ 6)$.

We remark that our result, especially the presence of ${\rm{M}}_{-4}$ and ${\rm{M}}_{-3}$, exhibits a connection between $L(\hf,\sym^2f)$ and quadratic fields. This connection is foreshadowed by the remarkable result of Zagier (see \cite[Theorem 1]{Zagier1977}) that a modular form $\Phi_s(z)$ of weight $k$ is constructed with Fourier coefficients being infinite linear combinations of zeta functions of quadratic fields such that
$$(\Phi_s,f)=\f{i^k\pi}{2^{k-3}(k-1)}\f{\Gamma(s+k-1)}{(4\pi)^{s+k-1}}L(s,\sym^2f).$$
It is worth pointing out that Kohnen$-$Sengupta \cite{KohnenSengupta2002} and Fomenko \cite{Fomenko2005} used the following consequence (see \cite[(2)]{KohnenSengupta2002}\footnote{In Eq.\,(2) of \cite{KohnenSengupta2002}, $2^{3k-3}$ should be replaced by $2^{3k-4}$, while $I_k(-4,0;\hf)$ by $2I_k(-4,0;\hf)$. }
%Such corrections are crucial for obtaining the correct contribution $\rm{M}_{1}$.}
or \cite[(17)]{Fomenko2005}) of the above identity to investigate the first moment in question:
\alig{\label{KSZ}
&\sum_{f\in\mcal{H}_{2k}}L\Big(\f{1}{2},\sym^2f\Big)\\
\nn&=\f{i^k2^{k-2}}{\sqrt{\pi}}\f{\Gamma(k)}{\Gamma(k-\hf)}
\bigg\{\lim_{s\rarrow\hf}\bigg[2(I_k(0,2;s)+I_k(0,-2;s)\zeta(2s-1)+\f{i^k\Gamma(s+k-1)\zeta(2s)}{2^{2s+k-3}\pi^{s-1}\Gamma(k)}\bigg]\\
\nn&\heq+\sum_{t\geq1}\bigg[I_k\Big(t^2-4,t;\hf\Big)+I_k\Big(t^2-4,-t;\hf\Big)\bigg]L\Big(\hf,t^2-4\Big)\bigg\},
}
We point out that in the above formula $L(\hf,t^2-4)=L(\hf,\chi_{t^2-4})$ for $t=1,\,2$, and refer the reader to \cite{Zagier1977} for the definitions of $L(s,\Delta)$ and $I_k(\Delta,t;s)$ where $\Delta=t^2-4m$ ($m=1$ for our case). Without explicitly stating an asymptotic formula, Kohnen and Sengupta isolated the secondary main term $\rm{M}_{-4}$ (the ``$t=2$'' term in (\ref{KSZ})). Fomenko obtained an asymptotic formula with error term of size $O(k^{-1/2})$. 
%Both Kohnen-Zagier and Fomenko were not aware of that $\rm{M}_{-3}$ (the ``$t=3$'' term in (\ref{KSZ})) can also contribute in a significant way. 
By comparing (\ref{KSZ}) and Theorem \ref{MainThm} we can see that
%\footnote{Footnote for myself: It is easy to check that the $t=2$ contribution matches $\rm{M}_{-4}$. But it needs more effort to verify that the $t=3$ contribution matches $\rm{M}_{-3}$ (already done but needs about a page to explain). Ongoing work: see if one can get an alternative proof of our result purely using Kohnen$-$Sengupta's approach based on Zagier's Theorem 1.}
the total contribution from all $t\geq3$ in (\ref{KSZ}) is $O_{\sss B}(k^{-B})$ for any fixed $B>0$ and sufficiently large $k$, which is difficult to obtain by direct computation.

Lastly, we include an example that also exhibits the connection mentioned in the last paragraph. Based on \cite{Zagier1977} Dummigan \cite{Dummigan2001} gives the following formula on (non-central) critical values of $L(s,\sym^2f)$: for fixed $k$ and odd $r$ with $3\leq r\leq k-1$
\aligs{
{\rm{w}}_f L(r,\sym^2f)=-\f{(2\pi)^{2r}}{4}\f{\Gamma(k-r)}{\Gamma(k+r-1)}\beta\\
}
with
\aligs{\beta=
\begin{cases}
c_{1}\zeta(1-2r)+c_{-4}L(1-r,\chi_{-4})+c_{-3}L(1-r,\chi_{-3}),&\quad 3\leq r<k-1,\vspace*{1ex}\\
(c_{1}+2k/B_k)\zeta(1-2r)+c_{-4}L(1-r,\chi_{-4})+c_{-3}L(1-r,\chi_{-3}),&\quad r=k-1,
\end{cases}
}
where $B_k$ is the $k$-th Bernoulli number,
$$
c_{1}=p_{k,r}(2,1)+p_{k,r}(-2,1),\ c_{-4}=p_{k,r}(0,1),\ \mbox{and }c_{-3}=p_{k,r}(1,1)+p_{k,r}(-1,1),
$$
and $p_{k,r}(t,m)$ denotes the coefficient of $x^{k-r-1}$ in $(1-tx+mx^2)^{-r}$. In view of the above formula and our result, for the first moment of critical values of symmetric square $L$-functions at $\hf$ and odd $r$ between $3$ and $k-1$ for large weight $k$, the quadratic fields of discriminant $1$ (degenerate), $-4$, and $-3$, play a dominant role over quadratic fields of other discriminants.
%Since forms $f\in\mcal{H}_k$ of weight $k=2m$ behave differently according to the parity of $m$ (odd $m$ implies $L(\hf,f)=0$
%$$
%\lim_{k\rarrow\infty}\f{\Gamma\!\el(\f{k-\hf}{2}\er)}{\Gamma\!\el(\f{k+\hf}{2}\er)}\cdot k^{\hf}=\sqrt{2},
%$$
%The symmetric square $L$-function
%$$
%L(s,\sym^2f)=\zeta(2s)\sum_{n\geq1}\f{\lambda_f(n^2)}{n^s}=\f{\zeta(2s)}{\zeta(s)}L(s,f\times f)
%$$

%It is well-known that $\lambda_f(n)$'s are real and that $L(s,f)$ admits an entire continuation to the whole $s$-plane and the complete $L$-function $\Lambda(s,f)=(2\pi)^{-s}\Gamma(s+k-\thf)L(s,f)$ satisfies the functional equation
%\eqs{
%\Lambda(s,f)=i^{2k}\Lambda(1-s,f).
%}
\subsection*{Acknowledgements}The author thanks Professor Wenzhi Luo for stimulating conversations and helpful comments. The author also thanks Professors Dorian Goldfeld, Roman Holowinsky, and Kannan Soundararajan for their interest in this work.
\section{Preparation}
In this section, we set notations and gather some preparatory results to be used later. We write $e(z)$ for $e^{2\pi i z}$, use $B$ for large real numbers, and reserve $p$ for prime numbers.% We adopt the conventions that these constants may take different values in different paragraphs, and that implied constants may depend on $\ep$ or $\delta$.
%For example, the implied constants in Proposition \ref{MoDerAsymp} depend on $\delta$ but the dependence is not shown for it is not important for our purpose.
\subsection{Functional equations}
We include here functional equations of zeta-functions and $L$-functions which will be used at various points.

The Riemann zeta-function $\zeta(s)$ satisfies the functional equation
\eq{\label{FE_zeta}
\zeta(s)=2(2\pi)^{s-1}\sin\Big(\f{\pi s}{2}\Big)\Gamma(1-s)\zeta(1-s)\quad(s\neq0,1).
}

The periodic zeta-function function $F(s,a)=\sum_{n\geq 1}e(na)n^{-s}$ ($\re(s)>1$, $0<a<1$) satisfies the functional equation
\eq{\label{FE_periodic}
F(s,a)=\f{\Gamma(1-s)}{(2\pi)^{1-s}}\el\{e\Big(\f{1-s}{4}\Big)\zeta(1-s,a)+e\Big(\f{s-1}{4}\Big)\zeta(1-s,1-a)\er\},
}
where $\zeta(s,a)=\sum_{n\geq0}(n+a)^{-s}$ is the Hurwitz zeta-function.

For a primitive Dirichlet character $\chi$ modulo $q$, the $L$-function $L(s,\chi)$ is entire. Let $\delta(\chi)=(1-\chi(-1))/2$. Then the complete $L$-function $\Lambda(s,\chi)=(q/\pi)^{\f{s+\delta(\chi)}{2}}\Gamma(\f{s+\delta(\chi)}{2})L(s,\chi)$ satisfies the functional equation
\eq{\label{FE_Dirichlet}
\Lambda(s,\chi)=\ep(\chi)\Lambda(1-s,\lbar{\chi})
}
where the root number $\ep(\chi)=i^{-\delta(\chi)}\tau(\chi)q^{-\hf}$. In particular, $\ep(\chi_{-4})=\ep(\chi_{-3})=1$.

For the symmetric square $L$-function $L(s,\sym^2 f)$ let
\aligs{
L_\infty(s)=\pi^{-\f{3}{2}s}\Gamma\Big(\f{s+1}{2}\Big)\Gamma\Big(\f{s+k-1}{2}\Big)\Gamma\Big(\f{s+k}{2}\Big).
}
Then the complete $L$-function $\Lambda(s,\sym^2f)=L_\infty(s)L(s,\sym^2 f)$ satisfies the functional equation
\eq{\label{FE_SymSq}
\Lambda(s,\sym^2f)=\Lambda(1-s,\sym^2f).
}
%Since the complete $L$-function $\Lambda(s,\chi_{-4})=(\f{4}{\pi})^{\f{1+s}{2}}\Gamma(\f{s+1}{2})L(s,\chi_{-4})$ satisfies the functional equation
%$$
%\Lambda(s,\chi_{-4})=\Lambda(1-s,\chi_{-4}),
%$$
\subsection{Approximate functional equation}
Unlike previous studies we use
$$
L(s,\sym^2f)=\f{\zeta(2s)}{\zeta(s)}\sum_{n\geq1}\f{\lambda_f(n)^2}{n^s}
$$
to derive the following the approximate functional equation for $L(\thf,\sym^2f)$.
\lem{\label{AFE_SymmSq}Let $H(u)=e^{-u^4}u^{-1}$. We have
\eqs{
L\Big(\hf,\sym^2f\Big)=2\sum_{n\geq1}\f{\lambda_f(n)^2}{\sqrt{n}}V_k(n),
}
where
$$
V_k(y)=\f{1}{2\pi i}\int_{(3)}y^{-u}H(u)R(u)\f{\zeta(1+2u)}{\zeta(\hf+u)}\dd u
$$
and
$$R(u)=\f{L_\infty(\hf+u)}{L_\infty(\thf)}=\pi^{-\f{3}{2}u}R_a(u)=(2\pi^{\f{3}{2}})^{-u}R_b(u)$$
with
\aligs{
R_a(u)=\f{\Gamma\el(\f{k-\hf+u}{2}\er)}{\Gamma\el(\f{k-\hf}{2}\er)}\f{\Gamma\el(\f{k+\hf+u}{2}\er)}{\Gamma\el(\f{k+\hf}{2}\er)}\f{\Gamma(\f{3}{4}+\f{u}{2})}{\Gamma(\f{3}{4})}\qand R_b(u)=\f{\Gamma(k-\hf+u)}{\Gamma(k-\hf)}\f{\Gamma(\f{3}{4}+\f{u}{2})}{\Gamma(\f{3}{4})}.
}
}
The proof is standard (see \cite[Theorem 5.3]{IwaniecKowalski2004}). We will use in the sequel whichever version of $R(u)$ that is more convenient for computation.
\subsection{Petersson's formula}For a nonempty Hecke eigenbasis $\mcal{H}_k$, Petersson's formula (see \cite{IwaniecKowalski2004}) states that
% (see \cite[Corollary 14.23 ]{IwaniecKowalski2004}, for example)
$$
\f{\Gamma(k-1)}{(4\pi)^{k-1}}\!\sum_{f\in\mcal{H}_k}\!\f{\lambda_f(m)\lambda_f(n)}{(f,f)}=\delta_{m,n}+2\pi i^{-k}\sum_{c\geq1}\f{S(m,n;c)}{c}J_{k-1}\!\el(\f{4\pi\sqrt{mn}}{c}\er)\!.
$$
Here $\delta_{m,n}$ is the Kronecker delta symbol; $S(m,n;c)$ denotes the classical Kloosterman sum
$$
S(m,n;c)=\!\sum_{\subs{x(\mmod\,c)\\x\lbar{x}\equiv 1 (\mmod\,c)}}\!e\!\el(\f{mx+n\lbar{x}}{c}\er);
$$
%which satisfies Weil's bound $|S(m,n;c)|\leq \tau(c)(m,n,c)^{\hf}c^\hf$;
$J_\nu(z)$ is the Bessel function of the first kind of order $\nu$. For convenience, we will denote the summation symbol in $S(m,n;c)$ by ${\sum'_{x(c)}}$.
%
%
%
%\subsection{Some properties of $J_{\nu}(x)$}
\subsection{An integral representation of $J_{\nu}(x)$}
%For $\nu\geq0$, we have
%\eqs{
%|J_\nu(x)|\leq1\quad\mbox{for }\ x\in\mbb{R}%\!\quad\mbox{(see \cite[(10.14.1)]{Ol})}
%}
%and by \cite[3.31(1)]{Watson1996} and \cite[(5.6.1)]{Ol} that
%\eqs{
%|J_\nu(x)|\leq \f{e}{\sqrt{2\pi}}\f{1}{\sqrt{\nu+1}}\!\el(\f{ex}{2\nu+2}\er)^\nu.%^\nu\!\quad\mbox{(by \cite[3.31(1)]{Watson1996} and \cite[(5.6.1)]{Ol}})
%}
%Combining these, we get the following bound which is sufficient for our purpose: for any positive integer $k$ it holds that
%\eq{\label{Jbound2}
%J_{2k-1}(x)\ll\begin{cases}2^{-k}k^{-\f{3}{2}}x, &\mbox{ if }0<x<k,\vspace*{0.5ex}\\ 1, &\mbox{ if }x\geq k. \end{cases}
%}
We also need the following Mellin-Barnes integral
\eq{\label{MBJ}
J_\nu(x)=\f{1}{2\pi i}\int_{\sss (\sigma)}x^{-s-1}2^s\f{\Gamma(\f{\nu+1+s}{2})}{\Gamma(\f{\nu+1-s}{2})}\dd s
}
for $x>0$ and $-1-\re(\nu)<\sigma<0$ (see \cite[p.\,82]{MOS1966} or \cite[p. 226]{Luo2015}).
\subsection{A complex integral}
In the analysis of the contribution from off-diagonal terms resulting from Petersson's formula, we will encounter the following complex integral
\eqs{
\mcal{I}(u,y)=
\f{1}{2\pi i}\ints{(3)}\f{\Gamma\!\el(\f{k-\hf-u-z}{2}\er)}{\Gamma\!\el(\f{k+\hf+u+z}{2}\er)}\Gamma(z)
\cos\!\Big(\f{\pi z}{2}\Big)\,y^{-z}\dd z.
}
\lem{\label{CXintegral}For $y>0$ and $\hf<\re(u)<k-4$,
\eqs{
\mcal{I}(u,y)=
\begin{cases}
\ds\f{\Gamma\!\el(\f{k-\hf-u}{2}\er)}{\Gamma\!\el(\f{k+\hf+u}{2}\er)} F\bigg(\f{k-\thf-u}{2},\f{-k+\tf{3}{2}-u}{2};\hf;\f{y^2}{4}\bigg),&y<2,\vspace*{1.5ex}\\
\ds\f{i^k 2^{2u}\cos(\f{\pi}{2}(\hf+u))}{\sqrt{\pi}}\Gamma(u)\f{\Gamma(k-\hf-u)}{\Gamma(k-\hf+u)},&y=2\vspace*{1.5ex},\\
\ds\f{2i^k\cos(\f{\pi}{2}(\hf+u))}{y^{k-\hf-u}}\f{\Gamma(k-\hf-u)}{\Gamma(k)}\ds F\bigg(\f{k-\thf-u}{2},\f{k+\thf-u}{2};k;\f{4}{y^2}\bigg),&y>2,
\end{cases}
}
where $F(a,b;c;z)$ denotes the Gauss hypergeometric function.
}
\pf{Since $\mcal{I}(u,y)$ clearly is holomorphic in $u$ for $\hf<\re(u)<k-4$, we only need to establish the result for $u>1$, say.
First we have the Mellin inversion (see \cite[Vol.\,1, (2.5.3.10)]{PBM1986} or \cite[(1.3.1)]{Erdelyi1954})
$$
\cos(x)=\f{1}{2\pi i}\int_{({\sss \hf})}\Gamma(z)\cos\!\Big(\f{\pi z}{2}\Big)\,x^{-z}\dd z.
$$
We also note that
$$
\f{\Gamma(\f{k+s}{2})}{\Gamma(\f{k-s}{2})}=\int_0^\infty J_{k-1}(x)\Big(\f{x}{2}\Big)^s\dd x,\quad -k<\re(s)<-\thf.
$$
Shifting the contour of $\mcal{I}(u,y)$ to $\re(z)=\thf$, we get
\aligs{
\mcal{I}(u,y)&=\lim_{T\rarrow\infty}\int_{\sss \hf-iT}^{\sss \hf+iT}\el(\int_0^\infty J_{k-1}(x)\Big(\f{x}{2}\Big)^{-\hf-u-z}\dd x\er)\Gamma(z)\cos\!\Big(\f{\pi z}{2}\Big)\,y^{-z}\dd z\\
&=2^{\hf+u}\int_0^\infty J_{k-1}(x)x^{-\hf-u}\el(\f{1}{2\pi i}\int_{({\sss \hf})}\Gamma(z)\cos\!\Big(\f{\pi z}{2}\Big)\Big(\f{xy}{2}\Big)^{-z}\dd z\er)\dd x\\
&=2^{\hf+u}\int_0^\infty J_{k-1}(x)\cos\!\Big(\f{xy}{2}\Big)\,x^{-\hf-u}\dd x.
}
Then the lemma follows from \cite[Vol.\,2, (2,12.15.3), (2.12.15.4), and (2.12.15.19)]{PBM1986} (or see \cite[(1.12.13) and (6.8.11)]{Erdelyi1954}).
}
\section{The diagonal contribution}
In this section we set up the proof of Theorem \ref{MainThm} and analyze the diagonal contribution.

By Lemma \ref{AFE_SymmSq} and Petersson's formula, we have
\aligs{
\sum_{f\in\mcal{H}_k}{\rm{w}}_fL\Big(\hf,\sym^2f\Big)
&=2\sum_{n\geq1}\f{V_k(n)}{\sqrt{n}}\sum_{f\in\mcal{H}_k}{\rm{w}}_f\lambda_f(n)^2\\
&=2\sum_{n\geq1}\f{V_k(n)}{\sqrt{n}}\bigg\{1+2\pi i^{-k}\sum_{c\geq1}\f{S(n,n;c)}{c}J_{k-1}\!\el(\f{4\pi n}{c}\er)\bigg\}\\
&=:\mcal{D}+\mcal{J},
}
where
$$
\mcal{D}=2\sum_{n\geq1}\f{V_k(n)}{\sqrt{n}}\qand  \mcal{J}=4\pi i^{-k}\sum_{n\geq1}\f{V_k(n)}{\sqrt{n}}\sum_{c\geq1}\f{S(n,n;c)}{c}J_{k-1}\!\el(\f{4\pi n}{c}\er).
$$
Using the definition of $V_k$ (in Lemma \ref{AFE_SymmSq}) and shifting the contour to $\re(s)=-B<0$, we have
\alig{\label{L_zeta}
\mcal{D}&=\f{2}{2\pi i}\ints{(3)}H(u)R(u)\f{\zeta(1+2u)}{\zeta(\thf+u)}\sum_{n\geq1}\f{1}{n^{\hf+u}}\dd u\\
\nn&=\f{2}{2\pi i}\ints{(3)}H(u)R(u)\zeta(1+2u)\dd u\\
\nn&=2\res{u=0}H(u)R(u)\zeta(1+2u)+\f{2}{2\pi i}\ints{(-B)}\cdots\\
\nn&=:\mcal{D}_R+\mcal{D}_I.
}
\textit{Note}: We adopt the convention that after a contour shifting the subscript ``$_R$'' indicates the residue part while the subscript ``$_I$'' the part from the resulting integral.

For the double pole $u=0$, we can easily compute
\aligs{
\mcal{D}_R&=2\lim_{u\rarrow 0}\f{\dd}{\dd u}u^2H(u)(2\pi^{\f{3}{2}})^{-u}R_b(u)\zeta(1+2u)\\
&=\psi\Big(k-\hf\Big)+2\gamma+\hf\psi\Big(\f{3}{4}\Big)-\log(2\pi^{\f{3}{2}}).
}
On the other hand we have
\aligs{
\mcal{D}_I&=\f{2}{2\pi i}\ints{(-B)}H(u)(2\pi^{\f{3}{2}})^{-u}R_b(u)\zeta(1+2u)\dd u\\
&=\f{2}{2\pi i}\ints{(-B)}H(u)(2\pi^{\f{3}{2}})^{-u}\f{\Gamma(k-\hf+u)}{\Gamma(k-\hf)}\f{\Gamma(\f{3}{4}+\f{u}{2})}{\Gamma(\f{3}{4})}\zeta(1+2u)\dd u\\
&\ll_B k^{-B},
}
in view of the classical bound of $\zeta$ and the bound
$$
\f{\Gamma(x+z)}{\Gamma(x)}\ll_{\sss\re(z)} x^{\re(z)}\quad\mbox{for }z\mbox{ with fixed }\re(z)\mbox{ and large }x>0.
$$
\textit{Note}: We will make free use of such bounds in the sequel.

Hence the diagonal contribution is
\eq{\label{D}
\mcal{D}=\psi\Big(k-\hf\Big)+2\gamma+\hf\psi\Big(\f{3}{4}\Big)-\log(2\pi^{\f{3}{2}})+O_{\sss B}(k^{-B}).
}
\section{The off-diagonal contribution}
In this section, we deal with the more delicate off-diagonal contribution $\mcal{J}$.
\subsection{Partition of $\mcal{J}$}For further analysis we first divide $\mcal{J}$ into several parts (see (\ref{Jparts})). To begin with, we have by (\ref{MBJ})
$$
J_{k-1}\!\el(\f{4\pi n}{c}\er)=\f{c}{4\pi n}\f{1}{2\pi i}\ints{(-1)}\el(\f{c}{2\pi n}\er)^s\f{\Gamma\big(\f{k+s}{2}\big)}{\Gamma\big(\f{k-s}{2}\big)}\dd s.
$$
Then by inserting the definition of $V_k$ and opening the Kloosterman sums, we have
\aligs{
\mcal{J}%&=i^{-k}\sum_{n\geq1}\f{V_k(n)}{\sqrt{n}}\sum_{c\geq1}\f{S(n,n;c)}{c}J_{k-1}\!\el(\f{4\pi n}{c}\er)\\
=\f{i^{-k}}{2\pi i}\ints{(3)}H(u)R(u)\f{\zeta(1+2u)}{\zeta(\hf+u)}
\sum_{c\geq1}\sumprime{x(c)}\f{1}{2\pi i}\ints{ (-1)}\el(\f{c}{2\pi}\er)^s\f{\Gamma\big(\f{k+s}{2}\big)}{\Gamma\big(\f{k-s}{2}\big)}\sum_{n\geq1}\f{e\!\el(n\f{x+\lbar{x}}{c}\er)}{n^{\f{3}{2}+s+u}}\dd s\dd u.
}
Upon letting
\alig{\label{rxc}
r(x,c)\equiv x+\lbar{x}\,(\mmod\ c)\quad\mbox{with}\quad 1\leq r(x,c)\leq c\qand r_{x,c}=\f{r(x,c)}{c}
}
we get
$$
\mcal{J}=\f{i^{-k}}{2\pi i}\ints{(3)}H(u)R(u)\f{\zeta(1+2u)}{\zeta(\hf+u)}
\sum_{c\geq1}\sumprime{x(c)}\f{1}{2\pi i}
\ints{ (-1)}\el(\f{c}{2\pi}\er)^s\f{\Gamma\big(\f{k+s}{2}\big)}{\Gamma\big(\f{k-s}{2}\big)}F\Big(\f{3}{2}+s+u,r_{x,c}\Big)\dd s\dd u,
$$
where $F(s,a)$ denotes the periodic zeta-function as in \S\,2.1. It is natural to write
$$
\mcal{J}=\mcal{J}_1+\mcal{J}_2
$$
according to $r_{x,c}=1$ (for which $F(s,r_{x,c})=\zeta(s)$) and $1/c\leq r_{x,c}<1$ (in which case $F(s,r_{x,c})$ is entire). For $\mcal{J}_1$, we shift the $s$-integral to $\re(s)=-B<0$, pick up a simple pole at $s=-\thf-u$, and get
$$\mcal{J}_1=\mcal{J}_{R}+\mcal{J}_{1,I}.$$
For $\mcal{J}_2$, we also shift its $s$-integral to $\re(s)=-B<0$ but get only
$$\mcal{J}_2=\mcal{J}_{2,I}, $$
since there is no residue part. Hence we arrive at the partition of $\mcal{J}$
\alig{\label{Jparts}
\mcal{J}=\mcal{J}_R+\mcal{J}_I
}
where $\mcal{J}_I=\mcal{J}_{1,I}+\mcal{J}_{2,I}$. In the sequel, we analyze $\mcal{J}_R$ (see (\ref{J_R_SymmSq})) and $\mcal{J}_I$ (see (\ref{J_I_SymmSq})) separately.
\subsection{Treatment of $\mcal{J}_R$}
Define
$$
N(c)=\sumprime{x(c)\atop r_{x,c}=1}1\qand L(s,N)=\sum_{c\geq1}\f{N(c)}{c^s}\quad(\re(s)>2).
$$
Then using
$$
\res{s=-{\sss\hf}-u}\el(\f{c}{2\pi}\er)^s\f{\Gamma\big(\f{k+s}{2}\big)}{\Gamma\big(\f{k-s}{2}\big)}\zeta\Big(\f{3}{2}+s+u\Big)
=\f{(2\pi)^{\hf+u}}{c^{\hf+u}}\f{\Gamma\!\el(\f{k-\hf-u}{2}\er)}{\Gamma\!\el(\f{k+\hf+u}{2}\er)}
$$
and $R(u)=\pi^{-\f{3}{2}u}R_a(u)$, we have
%$$
%\mcal{J}_{R}=\f{i^{-k}\pi^{\f{5}{4}}}{2\Gamma(\tf{3}{4})}\f{1}{2\pi i}
%\ints{(3)}H(u)\f{\Gamma\!\el(\f{k-\hf+u}{2}\er)}{\Gamma\!\el(\f{k-\hf}{2}\er)}\f{\Gamma\!\el(\f{k-\hf-u}{2}\er)}{\Gamma\!\el(\f{k+\hf}{2}\er)}
%\Big(\f{4}{\pi}\Big)^{\f{3}{4}+\f{u}{2}}\Gamma(\tf{3}{4}+\tf{u}{2})\f{\zeta(1+2u)}{\zeta(\thf+u)}L(\thf+u,N)\dd u.
%$$
\alig{\label{J_R_SymmSq}
\mcal{J}_{R}=\f{\pi^{\f{5}{4}}}{2}\f{i^{-k}}{2\pi i}
\ints{(3)}H(u)\f{\Gamma\!\el(\f{k-\hf+u}{2}\er)}{\Gamma\!\el(\f{k-\hf}{2}\er)}\f{\Gamma\!\el(\f{k-\hf-u}{2}\er)}{\Gamma\!\el(\f{k+\hf}{2}\er)}
\Big(\f{4}{\pi}\Big)^{\f{3}{4}+\f{u}{2}}\f{\Gamma(\tf{3}{4}+\tf{u}{2})}{\Gamma(\tf{3}{4})}\f{\zeta(1+2u)}{\zeta(\thf+u)}L\Big(\hf+u,N\Big)\dd u.
}

In order to evaluate $\mcal{J}_{R}$ we need the following result.
\lem{\label{L_N}
We have
$$\zeta(2s)L(s,N)=\zeta(s)L(s,\chi_{-4})=\zeta_{\mbb{Q}(i)}(s).$$
}
\pf{It suffices to consider $\re(s)>2$. It is clear from (\ref{rxc}) that
$$
N(c)=\#\big\{x\,(\mmod\ c)\mid (x,c)=1,\ x^2\equiv-1\,(\mmod\ c)\big\}.
$$
By the characterization for an element of $(\mbb{Z}/c\mbb{Z})^\times$ to be a quadratic residue mod $c$ (see \cite[Theorem 5.1]{LeVeque1996}, for example)
$$
N(c)=
\begin{cases}
2^{\omega(n)}R(-1,n),&\quad\mbox{if } c=2^an\mbox{ with }(2,n)=1,\ a=0\mbox{ or }1,\\
0,&\quad\mbox{otherwise},
\end{cases}
$$
where $\omega(n)$ is the number of distinct prime divisors of $n$ and
$$
R(a,n)=
\begin{cases}
1,&\quad\mbox{if }a\mbox{ is a quadratic residue mod }n,\\
0,&\quad\mbox{otherwise}.
\end{cases}
$$
Thus we have
$$
L(s,N)=\sum_{a=0}^1\sum_{{\rm odd}\,n\geq1}\f{2^{\omega(n)}R(-1,n)}{(2^an)^s}
=\el(1+\f{1}{2^s}\er)\sum_{n\geq1}\f{2^{\omega(n)}R(-1,n)P(n)}{n^s},
$$
where $P(n)=1$ if $n$ is odd and $=0$ if $n$ is even. It is easy to see that $2^{\omega(n)}R(-1,n)P(n)$ is multiplicative and that for any odd prime power $p^\alpha$ ($\alpha\geq1$)
$$
2^{\omega(p^\alpha)}R(-1,p^\alpha)P(p^\alpha)=2R(-1,p).
%=\begin{cases}
%2,&\quad\mbox{if }p\equiv1\,(\mmod\ 4),\\
%0,&\quad\mbox{if }p\equiv3\,(\mmod\ 4).
%\end{cases}
$$
Hence we have the Euler products
\aligs{
\f{L(s,N)}{\zeta(s)}&=\el(1-\f{1}{2^s}\er)\prod_{p>2}\el(1-\f{1}{p^s}\er)\\
&\heq\times\el(1+\f{1}{2^s}\er)\prod_{p>2}\el(1+\f{2R(-1,p)}{p^s}+\f{2R(-1,p)}{p^{2s}}+\cdots\er)\\
&=\el(1-\f{1}{2^{2s}}\er)\prod_{p>2}\el(1+\f{2R(-1,p)-1}{p^s}\er)\\
&=\el(1-\f{1}{2^{2s}}\er)\prod_{p>2}\el(1+\f{\chi_{-4}(p)}{p^s}\er)
}
and
\aligs{
\f{L(s,N)}{L(s,\chi_{-4})\zeta(s)}=\el(1-\f{1}{2^{2s}}\er)\prod_{p>2}\el(1-\f{\chi_{-4}(p)}{p^s}\er)\el(1+\f{\chi_{-4}(p)}{p^s}\er)=\f{1}{\zeta(2s)}.
}
%and\vspace*{1ex}\\
%$
%\hspace*{5em}\ds\f{L(s,N)}{L(s,\chi_{-4})\zeta(s)}=\el(1-\f{1}{2^{2s}}\er)\prod_{p>2}\el(1-\f{\chi_{-4}(p)}{p^s}\er)\el(1+\f{\chi_{-4}(p)}{p^s}\er)=\f{1}{\zeta(2s)}.
%$
}
Applying Lemma \ref{L_N} to (\ref{J_R_SymmSq}) we get
\alig{
\label{J_R}
\mcal{J}_{R}&=\f{i^{-k}\pi^{\f{5}{4}}}{2\Gamma(\tf{3}{4})}\f{1}{2\pi i}
\ints{(3)}H(u)\f{\Gamma\!\el(\f{k-\hf+u}{2}\er)}{\Gamma\!\el(\f{k-\hf}{2}\er)}\f{\Gamma\!\el(\f{k-\hf-u}{2}\er)}{\Gamma\!\el(\f{k+\hf}{2}\er)}
\Lambda\Big(\hf+u,\chi_{-4}\Big)\dd u\\
\nn&=\f{i^{-k}\pi^{\f{5}{4}}}{4\Gamma(\tf{3}{4})}\res{u=0}H(u)\f{\Gamma\!\el(\f{k-\hf+u}{2}\er)}{\Gamma\!\el(\f{k-\hf}{2}\er)}
\f{\Gamma\!\el(\f{k-\hf-u}{2}\er)}{\Gamma\!\el(\f{k+\hf}{2}\er)}
\Lambda\Big(\hf+u,\chi_{-4}\Big)\\
\nn&=i^{-k}\sqrt{\f{\pi}{2}}L\Big(\hf,\chi_{-4}\Big)\f{\Gamma\!\el(\f{k-\hf}{2}\er)}{\Gamma\!\el(\f{k+\hf}{2}\er)},
}
where we used that the integrand of the $u$-integral is an odd function in $u$ due to the functional equation (\ref{FE_Dirichlet}).
\subsection{Treatment of $\mcal{J}_I$}
We still need a few steps to get a reformulation (\ref{ReformJ_I}) of $\mcal{J}_{I}=\mcal{J}_{1,I}+\mcal{J}_{2,I}$. Applying the functional equation (\ref{FE_zeta}) to $\mcal{J}_{1,I}$ followed by the substitution $z=-\hf-s-u$ we have
\aligs{
\mcal{J}_{1,I}&=\f{i^{-k}}{2\pi i}\ints{(3)}H(u)R(u)\f{\zeta(1+2u)}{\zeta(\hf+u)}
\sum_{c\geq1}\sumprime{x(c)\atop r_{x,c}=1}\f{1}{2\pi i}\ints{(-B)}\el(\f{c}{2\pi}\er)^s\f{\Gamma\big(\f{k+s}{2}\big)}{\Gamma\big(\f{k-s}{2}\big)}\zeta\Big(\f{3}{2}+s+u\Big)\dd s\dd u\\
%&=\f{2i^{-k}}{2\pi i}\ints{(3)}H(u)(2\pi)^{\hf+u}R(u)\f{\zeta(1+2u)}{\zeta(\hf+u)}\sum_{c\geq1}\sumprime{x(c)\atop r_{x,c}=1}\\
%&\heq\times\f{1}{2\pi i}\ints{(-B)}c^s\f{\Gamma\big(\f{k+s}{2}\big)}{\Gamma\big(\f{k-s}{2}\big)}\Gamma\Big(-\hf-s-u\Big)\cos\!\Big(\f{\pi}{2}\Big(\hf+s+u\Big)\!\Big)\zeta\Big(-\hf-s-u\Big)\dd s\dd u\\
&=\f{2i^{-k}}{2\pi i}\ints{(3)}H(u)(2\pi)^{\hf+u}R(u)\f{\zeta(1+2u)}{\zeta(\hf+u)}\sum_{c\geq1}\f{1}{c^{\hf+u}}\sumprime{x(c)\atop r_{x,c}=1}\\
&\heq\times\f{1}{2\pi i}\int_{({\sss B-\f{7}{2}})}\f{1}{c^z}\f{\Gamma\!\el(\f{k-\hf-u-z}{2}\er)}{\Gamma\!\el(\f{k+\hf+u+z}{2}\er)}\Gamma(z)
\cos\!\Big(\f{\pi z}{2}\Big)\zeta(z)\dd z\dd u.
}
Similarly, by the functional equation (\ref{FE_periodic}) we get
\aligs{
\mcal{J}_{2,I}&=\f{i^{-k}}{2\pi i}\ints{(3)}H(u)(2\pi)^{\hf+u}R(u)\f{\zeta(1+2u)}{\zeta(\hf+u)}\sum_{c\geq1}\f{1}{c^{\hf+u}}\sumprime{x(c)\atop r_{x,c}<1}\\
&\heq\times\f{1}{2\pi i}\int_{({\sss B-\f{7}{2}})}\f{1}{c^z}\f{\Gamma\!\el(\f{k-\hf-u-z}{2}\er)}{\Gamma\!\el(\f{k+\hf+u+z}{2}\er)}\Gamma(z)
\el\{e\Big(\f{z}{4}\Big)\zeta(z,r_{x,c})+e\Big(\!-\f{z}{4}\Big)\zeta(z,1-r_{x,c})\er\}\dd z\dd u\\
&=\f{i^{-k}}{2\pi i}\ints{(3)}H(u)(2\pi)^{\hf+u}R(u)\f{\zeta(1+2u)}{\zeta(\hf+u)}\sum_{c\geq1}\f{1}{c^{\hf+u}}\sumprime{x(c)\atop r_{x,c}<1}\\
&\heq\times\f{1}{2\pi i}\int_{({\sss B-\f{7}{2}})}\f{1}{c^z}\f{\Gamma\!\el(\f{k-\hf-u-z}{2}\er)}{\Gamma\!\el(\f{k+\hf+u+z}{2}\er)}\Gamma(z)
\el\{e\Big(\f{z}{4}\Big)+e\Big(\!-\f{z}{4}\Big)\er\}\zeta(z,r_{x,c})\dd z\dd u,
}
since $1-r_{x,c}=r_{-x,c}$ for $x$ coprime with $c$ and if $x$ runs through reduced classes mod $c$ with $0<r_{x,c}<1$ then $-x$ also runs through reduced classes mod $c$ with $0<r_{-x,c}<1$.

Now we can rewrite $\mcal{J}_I=\mcal{J}_{1,I}+\mcal{J}_{2,I}$ as
\alig{\label{J_I_SymmSq}
\mcal{J}_I&=\f{2i^{-k}}{2\pi i}\ints{(3)}H(u)(2\pi)^{\hf+u}R(u)\f{\zeta(1+2u)}{\zeta(\hf+u)}\sum_{c\geq1}\f{1}{c^{\hf+u}}\sumprime{x(c)}\\
\nn&\heq\times\f{1}{2\pi i}\int_{({\sss B-\f{7}{2}})}\f{1}{c^z}\f{\Gamma\!\el(\f{k-\hf-u-z}{2}\er)}{\Gamma\!\el(\f{k+\hf+u+z}{2}\er)}\Gamma(z)
\cos\!\Big(\f{\pi z}{2}\Big)\zeta(z,r_{x,c})\dd z\dd u\\
\nn&=\f{2i^{-k}}{2\pi i}\ints{(3)}H(u)(2\pi)^{\hf+u}R(u)\f{\zeta(1+2u)}{\zeta(\hf+u)}\sum_{c\geq1}\f{1}{c^{\hf+u}}\sumprime{x(c)}\sum_{n\geq0}\\
\nn&\heq\times\f{1}{2\pi i}\int_{({\sss B-\f{7}{2}})}\f{\Gamma\!\el(\f{k-\hf-u-z}{2}\er)}{\Gamma\!\el(\f{k+\hf+u+z}{2}\er)}\Gamma(z)
\cos\!\Big(\f{\pi z}{2}\Big)(cn+r(x,c))^{-z}\dd z\dd u\\
\nn&=\f{2i^{-k}}{2\pi i}\ints{(3)}H(u)(2\pi)^{\hf+u}R(u)\f{\zeta(1+2u)}{\zeta(\hf+u)}\sum_{c\geq1}\f{1}{c^{\hf+u}}\sumprime{x(c)}\sum_{n\geq0}\mcal{I}(u;cn+r(x,c))\dd u,
}
by recalling the definition of $\mcal{I}(u;y)$ (see \S\,2.5). In view of Lemma \ref{CXintegral} it is natural to write
\alig{\label{ReformJ_I}
\mcal{J}_I=\mcal{J}_{I,\rm{i}}+\mcal{J}_{I,\rm{ii}}+\mcal{J}_{I,\rm{iii}},
}
where the summands correspond to the cases (i) $cn+r(x,c)<2$, (ii) $cn+r(x,c)=2$, and (iii) $cn+r(x,c)>2$, respectively.
\subsubsection{Treatment of $\mcal{J}_{I,\rm{i}}$}
Clearly $cn+r(x,c)=1$ if and only if $n=0$ and $r(x,c)=1$. Define
$$
M(c)=\sumprime{x(c)\atop r(x,c)=1}1\qand L(s,M)=\sum_{c\geq1}\f{M(c)}{c^s}\quad(\re(s)>2).
$$
Then
\aligs{
\mcal{J}_{I,\rm{i}}=\f{2i^{-k}}{2\pi i}\ints{(3)}H(u)(2\pi)^{\hf+u}R(u)\f{\zeta(1+2u)}{\zeta(\hf+u)}L\Big(\hf+u,M\Big)\,\mcal{I}(u;1)\dd u.
}

Similar to the treatment of $\mcal{J}_R$, we need the following result in order to evaluate $\mcal{J}_{I,\rm{i}}$.
\lem{\label{L_M}
We have
$$\zeta(2s)L(s,M)=\zeta(s)L(s,\chi_{-3})=\zeta_{\mbb{Q}(\sqrt{-3})}(s).$$
%where $\chi_{-3}(n)=(\f{-3}{n})$ is the Kronecker symbol.
}
\pf{
It is clear that
$$
M(c)=\#\big\{x\,(\mmod\ c)\mid (x,c)=1\mbox{ and }x^2-x+1\equiv0\,(\mmod\ c)\big\}.
$$
There are two steps to solve a quadratic equation $Ax^2+Bx+C\equiv0\,(\mmod\ D)$ with $(A,D)=1$: first solve $t^2\equiv B^2-4AC\,(\mmod\ 4D)$, $-D<t\leq D$, then solve $Ax\equiv (t-B)/2\,(\mmod\ D)$, $0\leq x<D$. For our case we first need to solve $t^2\equiv-3\,(\mmod\ 4c)$ for $-c<t\leq c$. It is clear that $t=c$ is a solution if and only if $c=1$ or $3$, for which $M(c)=1$. So we assume $c\neq1,\ 3$ and solve $t^2\equiv-3\,(\mmod\ 4c)$ for $-c<t<c$. Then it is not hard to see that $t=2c$ is not a solution and that
\alig{\label{half_soln}
\#\big\{\!-c<t<c\mid t^2\equiv-3\,(\mmod\ 4c)\big\}=\hf\,\#\big\{\!-2c<t<2c\mid t^2\equiv-3\,(\mmod\ 4c)\big\}.
}
%that if $0<t<c$ is a solution if and only if $c<2c-t<2c$ is a solution  that
%Also note that
%has a solution depends on whether $-3$ is a quadratic residue mod $4c$.
Alse note that for $n\geq1$ we have $(-3,3^n)=3$ and the equation $t^2\equiv-3\,(\mmod\ 3^n)$ has a solution if $n=1$ and no solution if $n>1$. Let $c=2^a3^bn$ with $(n,6)=1$. Then by virtue of (\ref{half_soln}) and \cite[Theorem 5.1]{LeVeque1996}, solving $t^2\equiv-3\,(\mmod\ 2^{a+2}3^bn)$ gives
$$
M(c)=
\begin{cases}
2^{\omega(n)}R(-3,n),&\quad\mbox{if } c=2^a3^bn\mbox{ with }(6,n)=1,\ a=b=0, \mbox{ or }a=0\mbox{ and }b=1,\\
0,&\quad\mbox{otherwise}.
\end{cases}
$$
Note that the above formula for $M(c)$ also includes the cases $c=1$ and $c=3$. Thus
$$
L(s,M)=\sum_{b=0}^1\sum_{n\geq1\atop (6,n)=1}\f{2^{\omega(n)}R(-3,n)}{(3^bn)^s}=\el(1+\f{1}{3^s}\er)\sum_{n\geq1}\f{2^{\omega(n)}R(-3,n)Q(n)}{n^s},
$$
where $Q(n)=\sum_{d\mid(6,n)}\mu(d)$. Again it is easy to check that $2^{\omega(n)}R(-3,n)Q(n)$ is multiplicative and for any odd prime power $p^\alpha$ ($p>3$ and $\alpha\geq1$)
$$
2^{\omega(p^\alpha)}R(-3,p^\alpha)Q(p^\alpha)=2R(-3,p).
%=\begin{cases}
%2,&\quad\mbox{if }p\equiv1\,(\mmod\ 4),\\
%0,&\quad\mbox{if }p\equiv3\,(\mmod\ 4).
%\end{cases}
$$
Hence the desired equality follows from an Euler product argument similar to that of Lemma \ref{L_N}.
%\aligs{
%\f{L(s,M)}{\zeta(s)}&=\el(1-\f{1}{2^s}\er)\el(1-\f{1}{3^{s}}\er)\prod_{p>3}\el(1-\f{1}{p^s}\er)\\
%&\heq\times\el(1+\f{1}{3^s}\er)\prod_{p>3}\el(1+\f{2R(-3,p)}{p^s}+\f{2R(-3,p)}{p^{2s}}+\cdots\er)\\
%&=\el(1-\f{1}{2^{s}}\er)\el(1-\f{1}{3^{2s}}\er)\prod_{p>3}\el(1+\f{2R(-3,p)-1}{p^s}\er)\\
%&=\el(1-\f{1}{2^{s}}\er)\el(1-\f{1}{3^{2s}}\er)\prod_{p>3}\el(1+\f{\chi_{-3}(p)}{p^s}\er)
%}
%\aligs{
%\f{L(s,M)}{L(s,\chi_{-3})\zeta(s)}=\el(1-\f{1}{2^{2s}}\er)\el(1-\f{1}{3^{2s}}\er)\prod_{p>3}\el(1-\f{\chi_{-3}(p)}{p^s}\er)\el(1+\f{\chi_{-3}(p)}{p^s}\er)=\f{1}{\zeta(2s)}.
%}
}
By Lemma \ref{CXintegral}, Lemma \ref{L_M}, that $R(u)=\pi^{-\f{3}{2}u} R_a(u)$, and the well-known linear transformation
$$
F(a,b;c;z)=(1-z)^{-a}F\Big(a,c-b;c;\f{z}{z-1}\Big),
$$
we have
\alig{\label{J_Ii1}
\mcal{J}_{I,\rm{i}}&=\f{2\sqrt{2\pi}i^{-k}}{\Gamma(\tf{3}{4})}\el(\f{\pi}{3}\er)^{\f{3}{4}}
\f{1}{2\pi i}\ints{(3)}H(u)\f{\Gamma\!\el(\f{k-\hf+u}{2}\er)}{\Gamma\!\el(\f{k-\hf}{2}\er)}\f{\Gamma\!\el(\f{k-\hf-u}{2}\er)}{\Gamma\!\el(\f{k+\hf}{2}\er)}
\Lambda\Big(\hf+u,\chi_{-3}\Big)\\
\nn&\hspace*{9em}\times\Big(\f{2}{\sqrt{3}}\Big)^u F\bigg(\f{k-\hf-u}{2},\f{-k+\f{3}{2}-u}{2};\hf;\f{1}{4}\bigg)\dd u\\
\nn&=\f{2\sqrt{2\pi}i^{-k}}{\Gamma(\tf{3}{4})}\el(\f{\pi}{3}\er)^{\f{3}{4}}
\f{1}{2\pi i}\ints{(3)}H(u)\f{\Gamma\!\el(\f{k-\hf+u}{2}\er)}{\Gamma\!\el(\f{k-\hf}{2}\er)}\f{\Gamma\!\el(\f{k-\hf-u}{2}\er)}{\Gamma\!\el(\f{k+\hf}{2}\er)}
\Lambda\Big(\hf+u,\chi_{-3}\Big)\\
\nn&\hspace*{9em}\times \Big(\f{2}{\sqrt{3}}\Big)^{k-\hf} F\bigg(\f{k-\hf-u}{2},\f{k-\hf+u}{2};\hf;-\f{1}{3}\bigg)\dd u
}
which becomes, since the integrand is odd in $u$ by the functional equation (\ref{FE_Dirichlet}),
\aligs{
&=\f{\sqrt{2\pi}i^{-k}}{\Gamma(\tf{3}{4})}\el(\f{\pi}{3}\er)^{\f{3}{4}}\Big(\f{2}{\sqrt{3}}\Big)^{k-\hf}
\res{u=0}\Bigg\{H(u)\f{\Gamma\!\el(\f{k-\hf+u}{2}\er)}{\Gamma\!\el(\f{k-\hf}{2}\er)}\f{\Gamma\!\el(\f{k-\hf-u}{2}\er)}{\Gamma\!\el(\f{k+\hf}{2}\er)}
\Lambda\Big(\hf+u,\chi_{-3}\Big)\\
&\hspace*{15em}\times F\bigg(\f{k-\hf-u}{2},\f{k-\hf+u}{2};\hf;-\f{1}{3}\bigg)\Bigg\}\\
&=\sqrt{2\pi}i^{-k}L\Big(\hf,\chi_{-3}\Big)
\Big(\f{2}{\sqrt{3}}\Big)^{k-\hf}F\bigg(\f{k-\hf}{2},\f{k-\hf}{2};\hf;-\f{1}{3}\bigg)
\f{\Gamma\!\el(\f{k-\hf}{2}\er)}{\Gamma\!\el(\f{k+\hf}{2}\er)}.
}
Next we analyize the behavior of $\big(\f{2}{\sqrt{3}}\big)^{k-\hf}
F\Big(\f{k-\hf}{2},\f{k-\hf}{2};\hf;-\f{1}{3}\Big)$. To this end we need \cite[Line 1 on p.\,54]{MOS1966}, that is, %for $0<x<\infty$
\aligs{
F\Big(a,b;\hf;-x\Big)&=\f{2^{a-b-1}}{\sqrt{\pi}}\Gamma\Big(a+\hf\Big)\Gamma(1-b)(1+x)^{-\f{a+b}{2}}\\
&\heq\times\bigg\{P^{b-a}_{a+b-1}\bigg(\sqrt{\f{x}{1+x}}\bigg)+P^{b-a}_{a+b-1}\bigg(\!-\sqrt{\f{x}{1+x}}\bigg)\bigg\}\quad (0<x<\infty)
}
as well as the asymptotic formula \cite[(3.9.1.2)]{Erdelyi1953}, that is,
\aligs{
P^\mu_\nu(\cos\theta)=\f{\Gamma(\nu+\mu+1)}{\Gamma(\nu+\f{3}{2})}\sqrt{\f{2}{\pi\sin\theta}}\,
\Big\{\cos\Big[\Big(\nu+\hf\Big)\theta-\f{\pi}{4}+\f{\mu\pi}{2}\Big]+O(\nu^{-1})\Big\},
}
where $\mu$ and $\nu$ are real, $0<\ep<\theta<\pi-\ep$, and $|\nu|\gg\ep^{-1}$. With the above formulas we see that
\aligs{
&\Big(\f{2}{\sqrt{3}}\Big)^{k-\hf}F\bigg(\f{k-\hf}{2},\f{k-\hf}{2};\hf;-\f{1}{3}\bigg)\\
&=\f{\Gamma\!\el(\f{k+\hf}{2}\er)\Gamma\!\el(1-\f{k-\hf}{2}\er)}{2\sqrt{\pi}}\Big\{P^0_{k-\f{3}{2}}\Big(\hf\Big)+P^0_{k-\f{3}{2}}\Big(\!-\hf\Big)\Big\}\\
&=-\Big(\f{2}{\sqrt{3}}\Big)^{\hf}i^k\f{\Gamma\big(k-\hf\big)}{\Gamma(k)}\f{\Gamma\!\el(\f{k+\hf}{2}\er)}{\Gamma\!\el(\f{k-\hf}{2}\er)}\el\{C(k)+O(k^{-1})\er\}
}
where
$$
C(k)=\cos\!\Big(\f{k\pi}{3}-\f{7\pi}{12}\Big)+\cos\!\Big(\f{2k\pi}{3}-\f{11\pi}{12}\Big)=
\begin{cases}
\sqrt{\f{3}{2}},&\quad\mbox{if }k\equiv2\,(\mmod\ 6),\vspace*{1ex}\\
0,&\quad\mbox{if }k\equiv4\,(\mmod\ 6),\vspace*{1ex}\\
-\sqrt{\f{3}{2}},&\quad\mbox{if }k\equiv0\,(\mmod\ 6).
\end{cases}
$$
Hence we have
\alig{\label{J_Ii2}
\mcal{J}_{I,\rm{i}}=3^{\f{1}{4}}\sqrt{2\pi}L\Big(\hf,\chi_{-3}\Big)\f{\Gamma(k-\hf)}{\Gamma(k)}\big[S(k)+O(k^{-1})\big],
}
where
$$
S(k)=\begin{cases}
-1,&\quad\mbox{if }k\equiv2\,(\mmod\ 6),\vspace*{1ex}\\
0,&\quad\mbox{if }k\equiv4\,(\mmod\ 6),\vspace*{1ex}\\
1,&\quad\mbox{if }k\equiv0\,(\mmod\ 6).
\end{cases}
$$
\subsubsection{Treatment of $\mcal{J}_{I,\rm{ii}}$}
Clearly $cn+r(x,c)=2$ if and only if $c=n=1$ and $r(x,c)=1$ or $n=0$ and $r(x,c)=2$ for all $c\geq2$. By Lemma \ref{CXintegral} and that $R(u)=(2\pi^\f{3}{2})^{-u}R_b(u)$, we have
\alig{\label{J_Iii}
\mcal{J}_{I,\rm{ii}}&=\f{2i^{-k}}{2\pi i}\ints{(3)}H(u)(2\pi)^{\hf+u}(2\pi^{\f{3}{2}})^{-u}\f{\Gamma(k-\hf+u)}{\Gamma(k-\hf)}\f{\Gamma(\f{3}{4}+\f{u}{2})}{\Gamma(\f{3}{4})}\f{\zeta(1+2u)}{\zeta(\hf+u)}\\
\nn&\heq\times\f{i^k2^{2u}\cos(\tf{\pi}{2}(\thf+u))}{\sqrt{\pi}}\f{\Gamma(k-\hf-u)}{\Gamma(k-\hf+u)}\Gamma(u)
\bigg\{1+\sum_{c\geq2}\f{1}{c^{\hf+u}}\sumprime{x(c),r(x,c)=2}1\bigg\}\dd u\\
\nn&=\f{2\sqrt{2}}{2\pi i}\ints{(3)}H(u)\Big(\f{4}{\sqrt{\pi}}\Big)^u\f{\Gamma(k-\hf-u)}{\Gamma(k-\hf)}\f{\Gamma(\f{3}{4}+\f{u}{2})}{\Gamma(\f{3}{4})}\Gamma(u)
\cos\!\Big(\f{\pi}{2}\Big(\hf+u\Big)\!\Big)\f{\zeta(1+2u)}{\zeta(\hf+u)}\\
\nn&\heq\times\bigg\{1+\sum_{c\geq2}\f{1}{c^{\hf+u}}\sumprime{x(c),r(x,c)=2}1\bigg\}\dd u\\
\nn&\ll_{\sss B} k^{-B},
}
by shifting the integral to $\re(u)=B$ for large $B$.
\subsubsection{Treatment of $\mcal{J}_{I,\rm{iii}}$}
By Lemma \ref{CXintegral}, the well-known symmetry $F(b,a;c;z)=F(a,b;c;z)$, and the Euler integral representation for $F(a,b;c;z)$ with $\re(c)>\re(b)>0$ and $|\arg(1-z)|<\pi$, we have for $y>2$
$$
\mcal{I}(u,y)=\f{i^k}{\sqrt{\pi}}\f{\cos(\tf{\pi}{2}(\hf+u))}{(\f{y}{2})^{k-\hf-u}}\f{\Gamma\Big(\f{k+\hf-u}{2}\Big)}{\Gamma\Big(\f{k+\hf+u}{2}\Big)}
\int_0^1 t^{\f{k-\hf-u}{2}-1}(1-t)^{\f{k+\hf+u}{2}-1}(1-\tf{4}{y^2}t)^{-\f{k+\hf-u}{2}}\dd t
$$
where the $t$-integral equals
\alig{\label{int01}
&\int_0^1 t^{\f{k-\hf-u}{2}-1}(1-t)^{u-1}\bigg(\f{1-t}{1-\tf{4}{y^2}t}\bigg)^{\f{k+\hf-u}{2}}\dd t\\
\nn&\ll \int_0^1 t^{\f{k-\hf-\re(u)}{2}-1}(1-t)^{\re(u)-1}\dd t\\
\nn&=\f{\Gamma\Big(\f{k-\hf-\re(u)}{2}\Big)\Gamma(\re(u))}{\Gamma\Big(\f{k-\hf+\re(u)}{2}\Big)}\\
\nn&\ll_{\sss \re(u)}k^{-\re(u)}.
}
Hence by shifting the contour to $\re(u)=B$ and using $R(u)=\pi^{-\f{3}{2}u} R_a(u)$ and (\ref{int01}) we get
\alig{\label{J_Iiii}
\mcal{J}_{I,\rm{iii}}
&=\f{2i^{-k}}{2\pi i}\ints{(3)}H(u)(2\pi)^{\hf+u}R(u)\f{\zeta(1+2u)}{\zeta(\hf+u)}\sum_{c\geq1}\f{1}{c^{\hf+u}}\sumprime{x(c)}\sum_{n\geq0\atop cn+r(x,c)>2}\mcal{I}(u,cn+r(x,c))\dd u\\
\nn&=\f{2\sqrt{2}}{2\pi i}\ints{(B)}H(u)\Big(\f{2}{\sqrt{\pi}}\Big)^u
\f{\Gamma\Big(\f{k-\hf+u}{2}\Big)}{\Gamma\Big(\f{k-\hf}{2}\Big)}\f{\Gamma\Big(\f{k+\hf-u}{2}\Big)}{\Gamma\Big(\f{k+\hf}{2}\Big)}
\f{\Gamma(\f{3}{4}+\f{u}{2})}{\Gamma(\f{3}{4})}\cos\!\Big(\f{\pi}{2}\Big(\hf+u\Big)\!\Big)\f{\zeta(1+2u)}{\zeta(\hf+u)}\\
\nn&\heq\times\Bigg\{\sum_{c\geq1}\f{1}{c^{\hf+u}}\sumprime{x(c)}\sum_{n\geq0\atop cn+r(x,c)>2}\el(\f{cn+r(x,c)}{2}\er)^{-k+\hf+u}\\
\nn&\heq\heq\times\int_0^1 t^{\f{k-\hf-u}{2}-1}(1-t)^{\f{k+\hf+u}{2}-1}\el(1-\f{4t}{(cn+r(x,c))^2}\er)^{-\f{k+\hf-u}{2}}\dd t\Bigg\}\dd u\\
\nn&\ll_{\sss B} k^{-B}\sum_{c\geq1}\f{1}{c^{\hf+B}}\sumprime{x(c)}\sum_{n\geq0\atop cn+r(x,c)>2}\el(\f{cn+r(x,c)}{2}\er)^{-k+\hf+B}\\
\nn&\ll_{\sss B} k^{-B}\Bigg\{\sum_{c\geq1}\f{\phi(c)}{c^{\hf+B}}+\sum_{c\geq1}\f{\phi(c)}{c^{\hf+B}}\sum_{n\geq1}\el(cn\er)^{\f{-k+\hf+B}{2}}\Bigg\}\\
\nn&\ll_{\sss B} k^{-B},
}
where $\phi$ denotes the Euler totient function.

Finally, Theorem 1 follows from collecting the diagonal contribution (\ref{D}) and the off-diagonal contribution (\ref{J_R}), (\ref{J_Ii1}), (\ref{J_Iii}), and (\ref{J_Iiii}).

\end{document}